# A multiscale analysis of instability-induced failure mechanisms in fiber-reinforced composite structures via alternative modeling approaches

*Fabrizio Greco[1*], Lorenzo Leonetti[1], Paolo Lonetti[1], Raimondo Luciano[2], Andrea Pranno[1]*

[1]Department of Civil Engineering, University of Calabria, Rende, Italy

[2]Department of Engineering, Parthenope University of Naples, Naples, Italy

*Corresponding author. Department of Civil Engineering, University of Calabria, 87036, Rende, Cosenza, Italy. E-mail address: f.greco@unical.it (Fabrizio Greco).

**Abstract**

Multiscale techniques have been widely shown to potentially overcome the limitation of homogenization schemes in representing the microscopic failure mechanisms in heterogeneous media as well as their influence on their structural response at the macroscopic level. Such techniques allow the use of fully detailed models to be avoided, thus resulting in a notable decrease in the overall computational cost at fixed numerical accuracy compared to the so-called direct numerical simulations. In the present work, two different multiscale modeling approaches are presented for the analysis of microstructural instability-induced failure in locally periodic fiber-reinforced composite materials subjected to general loading conditions involving large deformations. The first approach, which is of a semi-concurrent kind, consists in the "on-the-fly" derivation of the macroscopic constitutive response of the composite structure together with its microscopic stability properties through a two-way computational homogenization scheme. The latter one is a novel hybrid hierarchical/concurrent multiscale approach relying on a two-level domain decomposition scheme used in conjunction with a nonlinear homogenization scheme performed at the preprocessing stage. Both multiscale approaches have been suitably validated through comparisons with reference direct numerical simulations, by which the ability of the latter approach in capturing boundary layer effects has been demonstrated.

# 1 Introduction

Polymer matrix composites (PMCs) have recently gained more and more popularity in the fields of applied chemistry, materials science and engineering, by virtue of their superior physical and mechanical features compared to conventional materials, especially in terms of stiffness-to-weight and strength-to-weight ratios. Moreover, their capability to undergo large deformations have promoted their use to manufacture highly stretchable polymer-based microstructured metamaterials (see, for instance, [1]), including supersensitive strain tensors for structural health monitoring devices. Generally speaking, such composite structures, whose reinforcing phase can be made of either fibers or particles with different sizes, must be explicitly designed to cope with extreme environmental conditions (such as corrosion, thermal ageing, static and dynamic multiaxial loadings, fatigue, etc.), which often expose them to severe stress and strain gradients, thus facilitating the occurrence of different kinds of failure, also at large deformation regimes. In particular, when subjected to high strain levels, these materials may exhibit instability phenomena at the microstructural scale (i.e. the scale of micro-constituents). Such phenomena are often regarded as a precursor of failure, since they may strongly interact with the microscopic damage processes, such as fiber/matrix or particle/matrix debonding, matrix cracking, fiber breakage, delamination between different plies, etc. [2–4]. Due to the complexity of the interactions between the above-described mechanical phenomena, occurring at multiple spatial scales, as well as of their effect on the resulting mechanical behavior at the macroscopic level, the accurate determination of the overall nonlinear response of such composite materials is still regarded as a challenging task, being also crucial to design advanced materials in innovative high-performance applications. Furthermore, the rapid growth in the use of composite materials in engineering applications for structural and mechanical systems has motivated the improvement of existing numerical theories for predicting the mechanical response of such materials, especially in the nonlinear framework [5–7]. Commonly, a composite material is idealized by a homogeneous continuum whose properties were extracted by using by using a suitable averaging technique and, thus, to give the macroscopic overall response in a homogenized sense, different homogenization techniques are typically adopted (e.g. [8–12]). A straightforward approach usually adopted to numerically model the heterogeneous materials is called Direct Numerical Simulation (or DNS in brief), based on the explicit discretization of the whole structure, meshing all the heterogeneities, which inevitably leads to a high computational effort. Numerous computational strategies have been developed in the past to reduce the computational cost of a direct numerical simulation, such as the approach proposed by Ladevèze et al. [13], related to the so-called LATIN method (a general iterative non-linear solver for time-dependent problems [14]), as well as various domain composition techniques and parallel computations. Another alternative is the well-established multilevel computational procedure called $FE^2$, proposed by Feyel and Chaboche [15], which is based on a multilevel finite element technique that simply correlates the microstructural response with the macroscopic behavior by means of the introduction of discretized representative volume elements attached to each integration point of the discretized homogenized macrostructure. By solving the nonlinear boundary value problems on the linked RVEs within a nested procedure, the macroscopic stresses can be obtained by direct averaging of the microscopic stresses (e.g. [16,17]), but it is well known that an acceptable computational effort is allowed only in the presence of an effective core processor parallelization implemented in the solver procedure, made possible by the fact that all the nested microscopic problems are decoupled to each other. In most cases, first-order schemes have been adopted for such homogenization techniques, whose validity requires for the microstructural length to be kept infinitely small compared to the macroscopic characteristic length, whereas a number of second-order constitutive or geometrical relations have been proposed to account for microstructural scales comparable with the characteristic size of the macroscopic domain (see, for instance, [18]). Moreover, with the aim of overcoming the limitations of classical first-order homogenization schemes, a coupled-volume multiscale approach has been proposed in [19], relying on the linking between the numerical parameter "macrolevel mesh size" and the model parameter "microlevel cell size", which is able to avoid any dependency on either the mesh size or the cell size. Furthermore, to circumvent the computational effort related to nested calculations in $FE^2$ methods, other approaches have been introduced for the homogenization of heterogeneous materials with time-dependent or nonlinear behaviors, for instance, recently a multiscale approach for the analysis of elastoplastic and viscoplastic composites employing a PieceWise Uniform TFA was proposed in [20]. Further advanced multiscale approaches, characterized by the hybrid multiscale modeling have been proposed with the aim to drastically reduce the computational effort while preserving a good numerical accuracy and effectiveness (see, for instance, [21–24]). In the framework of the multilevel analysis of heterogeneous materials, it worth noting that the instability phenomena can occur at both the micro- and macro-scales influencing each other. In [25] the influence of voids on the stability properties of a highly porous elastic material at finite levels of strain was investigated by using a consistent homogenization theory approach. This work highlighted that the matrix material used can be stable (in the sense that it remains elliptic at all deformations) but the resulting homogenized medium can exhibit a loss of ellipticity at adequately high strain levels. Moreover, the bifurcation of fiber-reinforced layered composite materials was investigated in [26] and the rigorous connection between loss of rank-one convexity of the homogenized tangent moduli and onset of bifurcation instability at the microscopic scale in finitely strained periodic elastomers was shown in [27]. In particular, it was shown that if the wavelength of the bifurcation eigenmode is much greater than the unit cell size, the corresponding instability of the periodic principal solution can be detected at the point of loss of ellipticity condition of the one-cell homogenized tangent moduli of the solid. Based on these results, in [28] the onset-of-failure surfaces in stress and strain spaces for periodic solids of infinite extent were defined, and subsequently applied to honeycombs ([29]), fiber-reinforced composites under combined normal and shear strains ([30]) and, more recently, to two-dimensional periodic elastomers with circular inclusions of variable stiffness-to-matrix ratios ([31]). The above investigations have pointed out that the onset of the macroscopic instability, characterized by wavelengths significantly larger than the microstructure characteristic size, can be predicted by the loss of ellipticity analysis requiring merely the evaluation of the one-cell homogenized tensor of elastic moduli, whereas the prediction of local instability modes with a finite wave-length is a more difficult task due to the theoretically infinite nature of the analysis domain, thus requiring sophisticated numerical techniques such as the Bloch wave stability analysis used in conjunction with the finite element method (see, for instance, [32–37]). Alternatively, the macroscopic constitutive stability measures can be adopted in order to obtain conservative prediction of the microscopic stability region [38]. In some recent studies [39–41], the combined effects of the interaction between various microscopic failure modes under large deformations (damage, fracture and instability, for instance) have been analyzed by also modeling unilateral self-contact between crack surfaces, and highlighting that a detailed microscopic continuum analysis with an appropriate modeling of different sources of nonlinearities is usually required, at the expense of a huge computational effort. Therefore, a consistent numerical modeling for the prediction of the compressive strength of composite materials must involve a multilevel analysis, by coupling a stability analysis at microscopic level with a macroscopic structural analysis. From this point of view, in order to avoid an explicit modeling of all microstructural details of the composite solid, multiscale techniques in conjunction with a first-order homogenization scheme can be adopted overcoming some limitations characterizing pure first-order homogenization approaches. A multilevel finite element technique, also called computational homogenization (FE$^2$), that considers two nested continuum models needing constitutive assumptions only at the local level, was implemented by Nezamabadi et al. in [42,43] proved a similar connection between instability at the local level and macroscopic loading noticed in previous referenced works. Further investigations can be found in [44], in which a numerical technique, combining the computational homogenization analysis and the asymptotic numerical method (ANM), is presented to deal with instability phenomena in the context of heterogeneous materials for which buckling may occur at either macroscopic or microscopic scales. Additional studies were undertaken by Miehe et al. [45], developing a general framework to examine, from both theoretical and computational points of view, the instability problems.

The previously mentioned multiscale techniques often suffer from some drawbacks, first of all their scarce predictive capabilities in the presence of strong coupling between micro- and macro-scales or of strong stress gradient associated with boundary layer effects, thus further investigation are still required to overcome such difficulties.

In the present work, two different multiscale modeling approaches are compared for the analysis of microstructural instability-induced failure in periodic fiber-reinforced composite materials in a large deformation context highlighting advantages and drawbacks of each implementation. The first approach, which is a semi-concurrent multiscale model, consists in the derivation of the macroscopic constitutive response of the composite structure together with their microscopic stability properties through a two-way computational homogenization scheme based on a coupled-volume scheme, similar to that proposed by Gitman et al. [19]. The second multiscale approach, introduced to overcome the limitations of semi-concurrent approaches (mainly, the inability to account for the boundary layer effects), is a hybrid hierarchical/concurrent multiscale approach relying on a two-level domain decomposition scheme implemented in conjunction with a nonlinear homogenization scheme performed at the preprocessing stage. Both multiscale approaches were firstly validated through comparisons with reference direct numerical simulations (DNS), based on the explicit modeling of the heterogeneities, especially in terms of instability critical load factors (obtained by means of a rigorous microscopic stability analysis). The excellent accuracy of the proposed hybrid multiscale model in predicting the failure onset in terms of loss of stability was proved together with its high capability in capturing boundary layer effects associated with complex loading and constraints. This work can be regarded as an extension of the previous work [21] by some of the authors, in which the main ingredient of the present hybrid multiscale model (i.e. the two-level domain decomposition scheme) was adopted to analyze the mechanical behavior of bioinspired nacre-like composites. In particular, a novel finite-element based numerical approach is here proposed, by coupling the previously mentioned two-level domain decomposition scheme with a rigorous microscopic stability analysis for the failure prediction of composite materials characterized by variegated microstructural geometries, including both continuous and discontinuous fiber reinforced ones.

## 2 Theoretical formulation of the microscopic stability problem for periodic composite solids

The theoretical formulation of the microscopic stability problem for perfectly periodic heterogeneous solids subjected to general loading conditions involving large deformations is here recalled, with reference to the special case of incrementally linear micro-constituents. In particular, Section 2.1 summarizes the main equations involved in the nonlinear homogenization problem, whereas Section 2.2 is devoted to the discussion of the microscopic stability condition.

Let us consider a perfectly periodic microstructured solid described by a unit cell attached to a generic material point $\bar{X}$ of the corresponding homogenized solid, as depicted in Fig. 1. The deformation of the considered microstructure is defined by the nonlinear mapping $x(X):V_{(i)} \to V$, relating points $X$ of the initial microstructural configuration $V_{(i)}$ to points $x$ of the current one, $V$. The deformation gradient at the material point $X$ is $F(X) = \partial x(X)/\partial X$, whereas the corresponding displacement field is $u(X) = x(X) - X$.

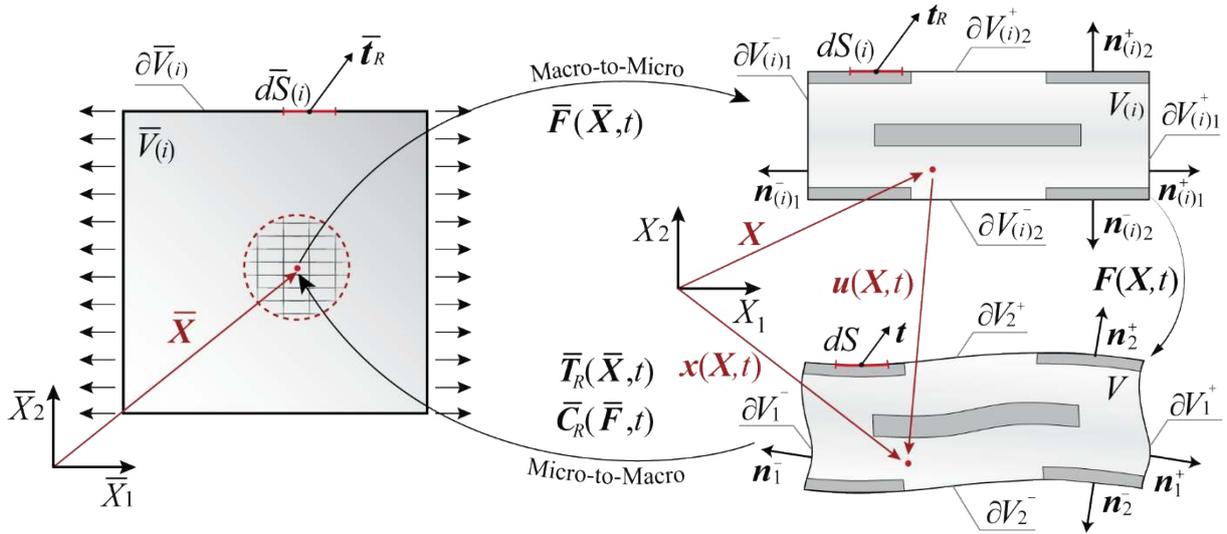

**Fig**.1 Homogenized solid and corresponding unit cells in the undeformed and deformed configurations attached to a generic material point $\bar{X}$.

### 2.1 Nonlinear homogenization problem: micro-to-macro coupling and macroscopic response

Each microstructural constituent is characterized by a rate independent material model, whose constitutive response at a microscopic point $X$ is described by an incrementally linear constitutive law:

$$\dot{T}_R = C^R(X,F)[\dot{F}], \tag{1}$$

relating the rate of deformation gradient tensor, $\dot{F}$, to the rate of the first Piola-Kirchhoff stress tensor, $\dot{T}_R$, via the nominal tangent moduli tensor $C^R$. Moreover, in the case of hyperelastic micro-constituents, the nominal stress tensor and the corresponding moduli tensor can be defined, respectively, as the first and second derivatives of the strain energy density function $W(X,F)$ with respect to $F$. As commonly assumed in the context of first-order homogenization schemes, the homogenized constitutive response of the considered microstructure relies on an equilibrium state which neglects body forces, resulting in a divergence-free local stress field $T_R$ in $V_{(i)}$.

In agreement with what reported in Section 3, the following expressions for the macroscopic nominal stress tensor $\bar{T}_R$ and the macroscopic deformation gradient tensor $\bar{F}$ are introduced:

$$\bar{T}_R = \frac{1}{|V_{(i)}|} \int_{\partial V_{(i)}} t_R(X) \otimes X \, dS_{(i)}, \quad \bar{F} = \frac{1}{|V_{(i)}|} \int_{\partial V_{(i)}} x(X) \otimes n_{(i)} \, dS_{(i)}, \tag{2}$$

with $\otimes$ denoting the tensor product, whereas $t_R$ and $n_{(i)}$ are the traction field and the outward normal, respectively, both evaluated at points $X$ belonging to the external boundary $\partial V_{(i)}$ of the unit cell. The microscopic deformation $x(X)$ can be additively split into a linear part, $\bar{F}X$, corresponding to a homogeneous deformation, and a correction part, $w(X)$, usually named as fluctuation field, being associated with nonhomogeneous deformations as follow:

$$x(X) = \bar{F}X + w(X). \tag{3}$$

By applying Eq. (2)$_2$, the following integral constraint turns to be required for the fluctuation field:

$$\int_{\partial V_{(i)}} w(X) \otimes n_{(i)} \, dS_{(i)} = 0. \tag{4}$$

According to the periodic nature of the composite microstructure, the constraint (4) can be satisfied by enforcing periodic fluctuations:

$$w(X^+) = w(X^-) \tag{5}$$

leading to periodic deformation and antiperiodic tractions on the unit cell boundary, i.e.

$$\begin{cases} u(X^+) - u(X^-) = (\bar{F}(\bar{X}) - I)(X^+ - X^-), \\ t_R(X^+) = -t_R(X^-), \end{cases} \quad \text{on } \partial V_{(i)}, \tag{6}$$

where $(X^+, X^-)$ is a couple of points belonging to the opposite sides of the unit cell boundary, denoted as $\partial V_{(i)}^+$ and $\partial V_{(i)}^-$ with outward normal unit vectors $n^+ = -n^-$ at two associated points $X^+ \in \partial V_{(i)}^+$ and $X^- \in \partial V_{(i)}^-$, respectively, obtained by two translations parallel to the directions of the periodicity vectors spanning the undeformed unit cell $V_{(i)}$. The incremental homogenized response of the given microstructure is determined considering a quasi-static loading path $\bar{F}(t)$ starting from its undeformed configuration $V_{(i)}$, where $t \geq 0$ represents a load parameter increasing monotonically with increasing prescribed macroscopic load and the given microstructure in the current configuration occupies the region $V$. The associated equilibrium solution at the given macro-deformation gradient can be obtained by solving the following variational problem defined over the unit cell:

$$\int_{V_{(i)}} T_R(X, \bar{F} + \nabla w) \cdot \nabla \delta w \, dV_{(i)} = 0, \quad \forall \delta w \in H^1(V_{(i)\#}), \tag{7}$$

where $H^1(V_{(i)\#})$ is the Sobolev space of vector-valued functions which are periodic over the unit cell $V_{(i)}$, the subscript # appended to a region denoting the assumed periodicity properties on its boundary. It is worth noting that arbitrary rigid body motions can be prevented by introducing additional constraints on the unknown fluctuation field. The equilibrium state obtained via the Euler-Lagrange equations associated with the variational formulation (7) is characterized by anti-periodic tractions on the external boundary $\partial V_{(i)}$. Moreover, if uniqueness of the equilibrium solution is assumed for each value of the loading parameter, the macroscopic loading path $\bar{F}(t)$ is named as principal solution path.

The incremental homogenized response of the microstructure can be derived after solving the following incremental equilibrium problem defined over the same unit cell subjected to an incremental change in the macroscopic deformation gradient, denoted as $\dot{\bar{F}}(t)$:

$$\int_{V_{(i)}} C^R(X, \bar{F})[\dot{\bar{F}} + \nabla \dot{w}] \cdot \nabla \delta \dot{w} \, dV_{(i)} = 0, \quad \forall \delta \dot{w} \in H^1(V_{(i)\#}), \tag{8}$$

where $\dot{w}$ is the unknown incremental fluctuation field induced by $\dot{\bar{F}}(t)$, and $C^R$ denotes the nominal moduli tensor, whose spatial distribution inside the unit cell is obtained after solving the microscopic problem (7). Additional displacement constraints are introduced into the variational problem (8) in order to exclude arbitrary incremental rigid body motions. The incremental equilibrium state obtained via the Euler-Lagrange equations associated with Eq. (8) is characterized by anti-periodic incremental tractions on $\partial V_{(i)}$. After solving the incremental problem (8) for the microstructure, its homogenized constitutive response is writable as $\dot{\bar{T}}_R = \bar{C}^R(\bar{F})[\dot{\bar{F}}]$, where the homogenized tangent moduli tensor $\bar{C}^R(\bar{F})$ can be expressed by the following relation, obtained by exploiting the fundamental identity $\dot{\bar{T}}_R = \overline{\dot{T}_R}$ (see [46] for details about its derivation):

$$\bar{C}^R_{ijkl}(\bar{F}) = \frac{1}{|V_{(i)}|} \int_{V_{(i)}} C^R_{ijmn}(X,\bar{F})[I^{kl}_{mn} + \nabla \dot{w}^{kl}_{mn}] \, dV_{(i)}, \qquad (9)$$

where $\dot{w}^{kl}$ is the incremental fluctuation field induced by unit components of the incremental macroscopic deformation gradient $\dot{\bar{F}} = I^{kl}$ (with $I^{kl}_{mn} = \delta_{mn}\delta_{kl}$).

It is important to recall here that Eq. (9), obtained via the unit cell-based homogenization, is strictly valid only if the microstructural equilibrium configuration is incrementally stable, otherwise the homogenized constitutive response must be computed considering a larger representative volume element (RVE), made of a possibly infinite number of unit cells, in order to capture all possible unstable deformation patterns.

**2.2 Microscopic stability conditions**

In the case of hyperelastic micro-constituents with nonconvex strain energy function $W(F)$, the homogenized strain energy function $\bar{W}(\bar{F})$ can be obtain by solving the following minimization problem:

$$\bar{W}(\bar{F}) = \inf_{k \in N} \left\{ \min_{w \in H^1(k^N V_{(i)}\#)} \left\{ \frac{1}{k^N |V_{(i)}|} \int_{k^N V_{(i)}} W(X, \bar{F} + \nabla w) \, dV_{(i)} \right\} \right\}, \qquad (10)$$

defined over all possible ensembles of $k^N = [0,k]^N$ unit cells, with $N = \{2,3\}$ and $k$ an arbitrary integer.

In the presence of eventual micro-buckling mechanisms, the solution of the minimization problem (10) allows to determine the size of the representative volume element associated with the fluctuation field capturing the minimizing buckling mode. However, the direct application of Eq. (10) may require a huge computational effort, connected to the need of investigating a full space at the microscopic scale. Therefore, a one-cell homogenization is preferable for numerical applications, obtained by taking $k=1$ in Eq. (10). In general, such a homogenization problem provides only an upper bound for the macroscopic strain energy, denoted as $\bar{W}^1(\bar{F})$, which becomes coinciding with the actual value $\bar{W}(\bar{F})$ only in the absence of micro-instabilities.

It can be shown that the region of validity for one-cell homogenization, i.e. the region of the macro-strain space for which $\bar{W}^1(\bar{F}) = \bar{W}(\bar{F})$, can be found in a rigorous manner by means of a microstructural stability analysis. For a given RVE subjected to a prescribed macro-deformation $\bar{F}$, such a microstructural stability condition relies on the positive definiteness of the following stability functional, written in a full Lagrangian setting:

$$S(\bar{F}, \dot{w}) = \int_{k^N V_{(i)}} C^R(X, \bar{F} + \nabla w)[\nabla \dot{w}(X)] \cdot \nabla \dot{w}(X) \, dV_{(i)}, \tag{11}$$

for all incremental fluctuations $\dot{w}$ satisfying the periodicity conditions on the RVE boundaries, with the additional constraint $\nabla \dot{w}(X) \neq \mathbf{0}$. It follows that the critical load parameter $t_c$ associated with the primary instability is detected when the minimum eigenvalue of $S(\bar{F}, \dot{w})$, taken over all admissible incremental fluctuations periodic on the ensemble of $k^N$ unit cells, first vanishes, namely:

$$\Lambda(t_c) = \inf_{k \in N} \left\{ \min_{\dot{w} \in H^1(k^N V_{(i)\#})} \left\{ \frac{S(\bar{F}(t_c), \dot{w})}{\int_{k^N V_{(i)}} \nabla \dot{w} \cdot \nabla \dot{w} \, dV_{(i)}} \right\} \right\} = 0. \tag{12}$$

At this load level, the initially stable and unique principal solution loses its uniqueness due to the existence of eigenmodes (i.e. non-trivial incremental periodic solutions to the homogeneous incremental equilibrium problem). Therefore, the microscopic stability region, defined as $t \,|\, \Lambda(\bar{F}(t)) > 0$ and characterized by an identical deformation for all the unit cells inside the considered RVE, necessarily coincides with the region of validity for one-cell homogenization.

Finally, it is worth recalling that a classical macroscopic stability analysis, based on the strong ellipticity condition for the homogenized moduli tensor, although being less computationally expensive in practical numerical applications, is not accurate in general, providing a non-conservative estimation of the primary microscopic instability load in the occurrence of instability modes of local kind (see [2,3,38,40,47–50] for additional details).

# 3 Description of the adopted multiscale approaches

In this section, in order to delineate which class of multiscale approach is more effective in the determination of the macroscopic behavior of advanced composite materials (fiber reinforced and nacre-like composites) in a large deformation context, more information about the computational implementation of the proposed multiscale approaches are given. Specifically, a semi-concurrent model called coupled-volume approach has been adopted to evaluate the microscopic instability critical load factor in composite materials reinforced with long fibers; later, a novel hybrid multiscale model has been proposed, firstly, to evaluate the microscopic instability in composite materials with staggered microstructure.

## 3.1 Coupled-volume semi-concurrent multiscale approach

Generally speaking, the concept behind the semi-concurrent models can be summarized in four steps, performed in an iterative manner:

(i) *Macro-domain computation*: the heterogeneous material is described as homogeneous with effective properties without the necessity to introduce constitutive assumptions;

(ii) *Downscaling:* once the macroscopic domain is discretized, a microscopic boundary value problem is linked to every Gauss quadrature point, and the information about the macroscopic gradient deformation field (input) is transferred to the microscopic domain in terms of displacement boundary conditions;

(iii) *Micro-domain computation:* at the microscopic level, each micro-constituent of the heterogeneous material is described by a constitutive assumption allowing to solve the $n$ boundary value problems linked to the $n$ Gauss points in different manners (FEM, Voronoi cell FEM, Fast Fourier Transforms, etc.)

(iv) *Upscaling:* the homogenization procedure is performed for each iteration needed to the solver to achieve the macroscopic convergence in terms of macroscopic stress tensor and macroscopic tangent moduli tensor, which results in the homogenized relations that has to be transferred from the microscopic to the macroscopic level.

The proposed semi-concurrent multiscale approach for the microstructural stability analysis of locally periodic composite materials relies on a coupled-volume multiscale model [19] according to which two nested equilibrium problems are solved at the same time (as sketched in Fig. 2). The approach leaves the concept that a finite microscopic cell size should be linked to an infinitely small macroscopic material point and adopts the concept that the macroscopic and the microscopic mesh sizes are uniquely linked following the rule that the macroscopic element size is equal to the microscopic cell size.

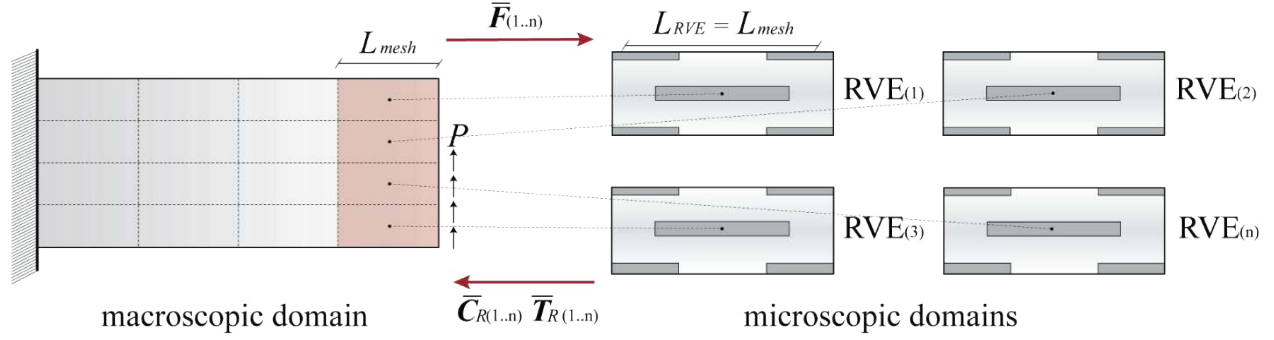

**Fig. 2** Schematic of the coupled-volume semi-concurrent multiscale approach.

This semi-concurrent approach is able to solve the mechanical problem associated with the macroscopic model in the absence of an explicitly defined macroscopic constitutive law, and it is also effective to investigate the behavior of materials subjected to strong nonlinearities for which the RVE cannot be found. Based on the information given above about the semi-concurrent multiscale modeling strategy, the implementation of the coupled-volume multiscale approach may be described by the steps which have been graphically summarized in Fig. 3.

The resulting two-scale equilibrium problem, written in an incremental form, reads as follows:

find $\dot{\bar{\boldsymbol{u}}} \in H^1(\bar{V}_{(i)}), \dot{\boldsymbol{w}} \in H^1(\bar{V}_{(i)} \times V_{(i)\#})$ such that:

$$\begin{aligned}
&\overbrace{\int_{\bar{V}_{(i)}} \dot{\bar{\boldsymbol{T}}}_R \cdot \nabla \delta \dot{\bar{\boldsymbol{u}}}\, d\bar{V}_{(i)} = \int_{\bar{V}_{(i)}} \dot{\bar{\boldsymbol{f}}} \cdot \delta \dot{\bar{\boldsymbol{u}}}\, d\bar{V}_{(i)} + \int_{\partial \bar{V}_{(i)}} \dot{\bar{\boldsymbol{t}}}_R \cdot \delta \dot{\bar{\boldsymbol{u}}}\, d\bar{S}_{(i)}}^{\text{Coarse-scale problem}} \quad \forall \delta \dot{\bar{\boldsymbol{u}}} \in H^1(\bar{V}_{(i)}) \\
&\underbrace{\int_{V_{(i)}} \boldsymbol{C}^R [\dot{\bar{\boldsymbol{F}}} + \nabla \dot{\boldsymbol{w}}] \cdot \nabla \delta \dot{\boldsymbol{w}}\, dV_{(i)} = 0 \qquad\qquad\qquad \forall \delta \dot{\boldsymbol{w}} \in H^1(V_{(i)\#})}_{\text{Fine-scale problem}}
\end{aligned} \tag{13}$$

where $H^1(\bar{V}_{(i)})$ is the Sobolev space of vector-valued functions defined on $\bar{V}_{(i)}$ and satisfying the homogeneous displacement-type conditions on its constrained portion of the boundary, whereas $H^1(\bar{V}_{(i)} \times V_{(i)\#})$ is the Sobolev space of vector-valued functions defined over the Cartesian product $\bar{V}_{(i)} \times V_{(i)\#}$ and satisfying the periodic fluctuation conditions over the unit cell boundaries.

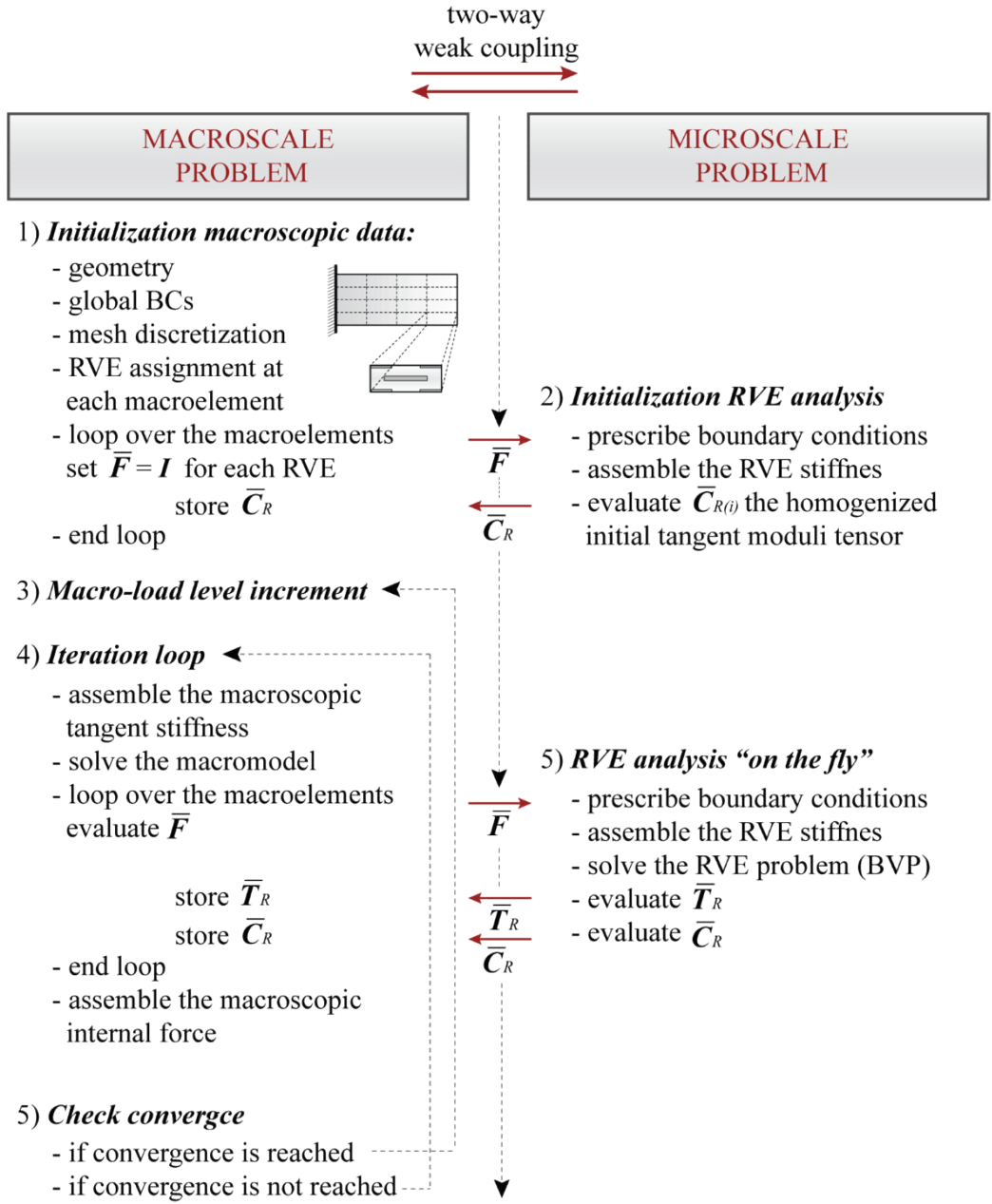

**Fig. 3** Implementation scheme of the coupled volume multiscale method.

Moreover, $\dot{\bar{f}}$ and $\dot{\bar{t}}_R$ denote, respectively, the body and surface incremental forces externally applied on the macroscopic body. The two-scale incremental equilibrium problem (13) turns to be two-way coupled. Indeed, the unit cell problem is driven by the incremental macroscopic deformation gradient $\dot{\bar{F}}$. Once the boundary value problem is solved for the unit cell, the rate of the macroscopic first Piola-Kirchhoff stress $\dot{\bar{T}}_R$ is computed via the first of Eqs.(2)$_1$ The information transfer between the two scales is sketched in Fig. 3.

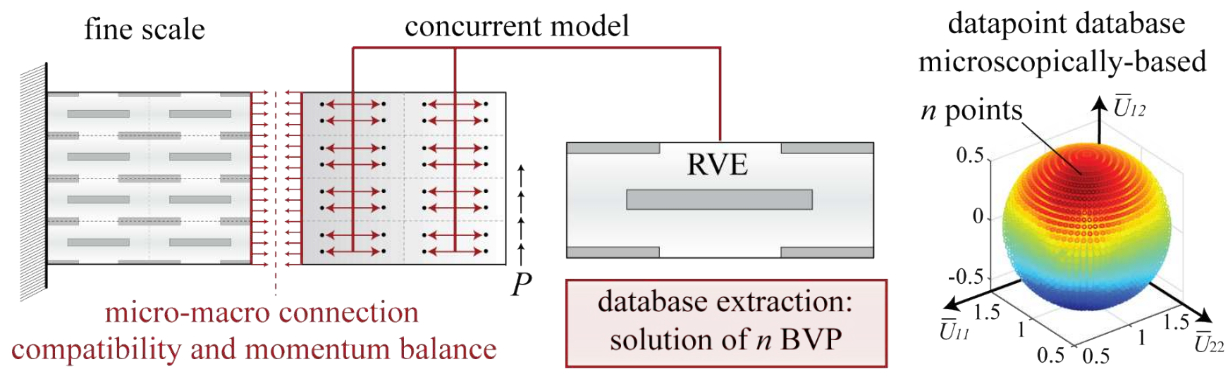

**Fig. 4** Schematic of the hybrid hierarchical/concurrent multiscale approach.

It is worth noting that the choice of using the deformation gradient $\boldsymbol{F}$ (together with its work-conjugate stress measure $\boldsymbol{T}_R$) rather than the Green-Lagrange strain $\boldsymbol{E}^{(2)}$ (and its work-conjugate stress measure $\boldsymbol{T}^{(2)}$) to drive the homogenization process is convenient for prescribing the periodic boundary conditions (BCs) on the RVE.

**3.2 Hybrid hierarchical/concurrent multiscale approach**

The key idea of the proposed hybrid hierarchical/concurrent multiscale approach is to combine a hierarchical and a concurrent approach overcoming the limitation of semi-concurrent approaches in capturing boundary layer effects, localization of deformation, coalescence of micro-cracks and other kind of nonlinearities. As shown in Fig. 4, the hybrid multiscale method can be implemented taking advantage of a concurrent approach identifying the so-called critical domains in which an explicit modeling of the microstructure (fine-scale domain) together with the associated nonlinear phenomena (material, geometrical, damage, friction, etc.) is employed, used in combination with a hierarchical approach implementing, in the remaining domain, a homogenized constitutive law in the form of nonlinear microscopically based database. Usually, the critical domains are characterized by a microstructural evolution involving notable geometrical or material nonlinearities, so that a numerical model able to describe all the different microscopic phenomena is required, leading to a huge computational effort; as a matter of fact, the identification of the regions where fine-scale computations are required represents the most important step of the entire procedure. Contrarily to the coupled-volume multiscale approach, in which a microscopic boundary value problem is solved for each RVE linked to a macro-element and for each time step of the macro-scale analysis, in the hybrid multiscale approach this step has been now replaced with the application of a microscopically based nonlinear database extracted in a preprocessing step. Thus, as explained in the following, the macroscopic stress and the macroscopic tangent moduli tensor are stored in a nonlinear database in the form of an unstructured point cloud and they are transferred to the macroscopic model by using an interpolation method.

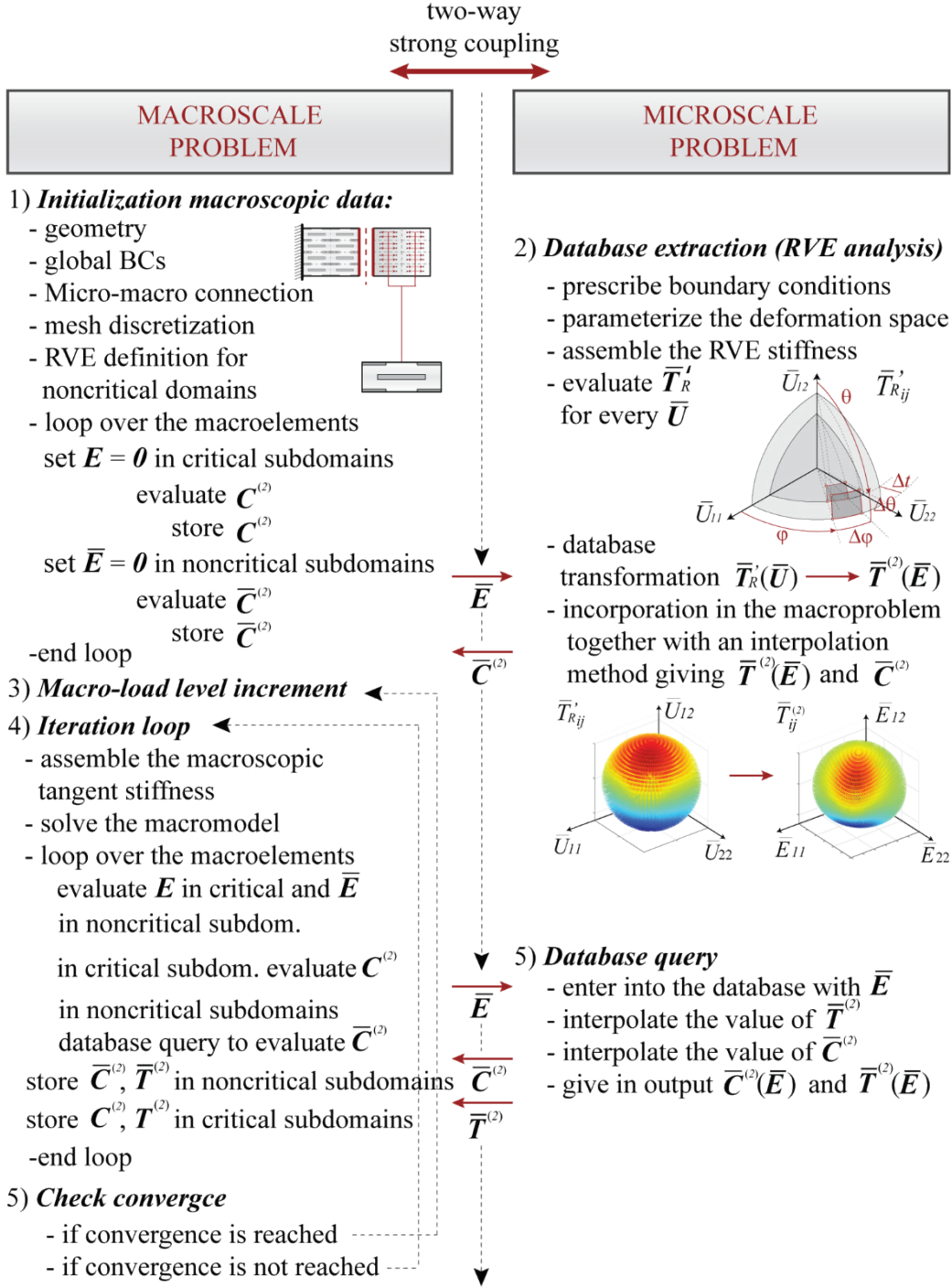

**Fig. 5** Implementation scheme of the hybrid hierarchical/concurrent multiscale method.

### 3.2.1 Numerical implementation

The procedure adopted to implement the hybrid multiscale approach, also sketched in Fig. 5, can be summarized in the following steps:

(i) *Macro-domain computation*: the heterogeneous domain is split in two subdomains: critical and noncritical. In critical subdomains every heterogeneity is modeled explicitly together with the nonlinear phenomena acting at the microscopic scale, assuming as known its constitutive properties. In noncritical subdomains no constitutive assumptions are required because the constitutive information is extracted from a nonlinear microscopically derived database. Once the macroscopic subdomains are discretized, compatibility and momentum balance are imposed on the common interface between critical and noncritical subdomains.

(ii) *Downscaling:* A representative volume element is identified to evaluate the homogenized behavior of the noncritical subdomains.

(iii) *Micro-domain computation and database extraction:* at the microscopic level, each micro-constituent of the representative volume element is equipped with a constitutive law. The RVE is analyzed to create a numerical data point in the macroscopic deformation space parametrized by using spherical coordinates and spaced by radial strain paths; the macroscopic stresses are evaluated for every time step and stored in a database matching the information on strains with information on stresses and tangent moduli.

(iv) *upscaling:* the extracted database is then incorporated in noncritical subdomains of the macroscopic model. During the procedure of solving the macroscopic problem, the microscopically based database is interrogated to extract the information in terms of macroscopic stress tensor and macroscopic tangent moduli tensor until the macroscopic convergence is reached.

More information about the step (iii) is given in the following section.

**3.2.2 Extraction of the microscopically derived database**

The microscopically based database represents a constitutive equation that receives in input a macroscopic strain tensor and gives in output a macroscopic stress tensor. Three alternative strategies can be adopted to extract the database:

i) the first one is based on the macroscopic First-Piola Kirchhoff tensor $\bar{T}_R$ defined as a function of the macroscopic gradient deformation tensor $\bar{F}$:

$$\bar{T}_R = \bar{T}_R\left(\bar{F}(t)\right). \tag{14}$$

This strategy leads to a simple definition of the RVE boundary conditions, but in a planar setting, $\bar{F}$ is a tensor whose components are: $\bar{F}_{11}$, $\bar{F}_{12}$, $\bar{F}_{21}$ and $\bar{F}_{22}$; for this reason, a numerical data point database, based on the macroscopic gradient tensor, must be parametrized on a four-dimensional space leading to a high computationally effort and leading also to the adoption of complex interpolation methods;

ii) the second one is based on the macroscopic Second-Piola Kirchhoff tensor $\bar{\boldsymbol{T}}^{(2)}$ defined as a function of its conjugate strain so-called Green-Lagrange strain tensor $\bar{\boldsymbol{E}}^{(2)}$:

$$\bar{\boldsymbol{T}}^{(2)} = \bar{\boldsymbol{T}}^{(2)}\left(\bar{\boldsymbol{E}}^{(2)}(t)\right). \tag{15}$$

This strategy leads to a complex definition of the RVE boundary conditions which are classically written in terms of the macroscopic gradient deformation tensor. In a planar setting, $\bar{\boldsymbol{E}}^{(2)}$ is a symmetric tensor whose components are: $\bar{E}_{11}$, $\bar{E}_{12} = \bar{E}_{21}$, $\bar{E}_{22}$ and, as a consequence, a numerical data point based on the macroscopic Green-Lagrange strain tensor can be parametrized on a three-dimensional space leading to a lower computational effort compared to the first strategy. However, the complexity involved in the imposition of boundary conditions makes this strategy unsuitable.

iii) the third strategy is able to take the advantages of the first two proposed strategies and to overcome their difficulties. This is based on the restricted First-Piola Kirchhoff stress tensor $\bar{\boldsymbol{T}}_R^{'}$ defined as a function of the macroscopic right stretch tensor $\bar{\boldsymbol{U}}$ and it leads to the extraction of the so-called reduced database:

$$\bar{\boldsymbol{T}}_R^{'} = \bar{\boldsymbol{T}}_R^{'}\left(\bar{\boldsymbol{U}}(t)\right). \tag{16}$$

This procedure exploits the following objectivity property [51] of the macroscopic strain energy density function $\bar{W}$:

$$\bar{W}(\bar{\boldsymbol{Q}}\bar{\boldsymbol{F}}) = \bar{W}(\bar{\boldsymbol{F}}) \quad \forall \bar{\boldsymbol{Q}} \in Orth^+, \tag{17}$$

where $\bar{\boldsymbol{Q}}$ is an arbitrary proper orthogonal tensor and $\bar{W}$ is defined as the unweighted volume average of the microscopic strain energy density function $W$:

$$\bar{W} = \frac{1}{|V_{(i)}|} \int_{V_{(i)}} W(\boldsymbol{X}, \boldsymbol{F}) dV. \tag{18}$$

By means of Eq. (17), which states that the rotational part of the macroscopic deformation $\bar{\boldsymbol{Q}}$ does not influence the constitutive properties of the homogenized material, one may rewrite the constitutive law in the following reduced form:

$$\bar{\boldsymbol{T}}_R(\bar{\boldsymbol{F}}) = \bar{\boldsymbol{R}}\bar{\boldsymbol{T}}_R^{'}(\bar{\boldsymbol{U}}), \tag{19}$$

where $\bar{T}_R'$ is a positive-definite symmetric tensor and $\bar{R}$ is the macroscopic rotation tensor obtained from the polar decomposition $\bar{F} = \bar{R}\bar{U}$ of the macroscopic gradient deformation tensor. In the light of the above considerations, throughout the database extraction phase one can assume $\bar{R} = I$ and, as a consequence, $\bar{F} = \bar{U}$, thus the reduced database can be parametrized on a three-dimensional space because, with reference to planar problems, $\bar{U}$ is a symmetric tensor whose components are $\bar{U}_{11}$, $\bar{U}_{12} = \bar{U}_{21}$, $\bar{U}_{22}$.

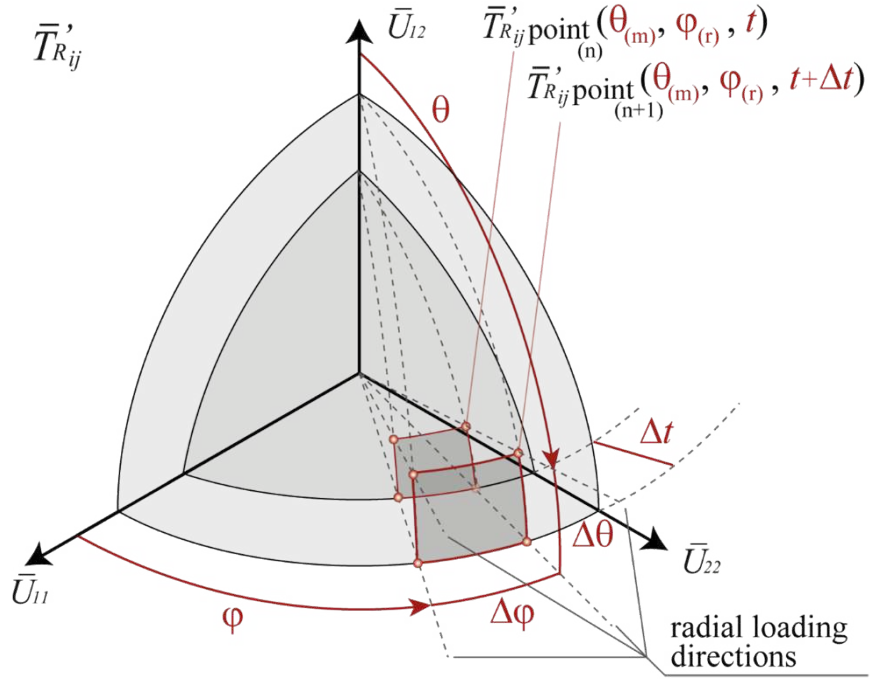

**Fig. 6** $\bar{U}_{11}$-$\bar{U}_{12}$-$\bar{U}_{22}$ strain space parametrized by spherical coordinates and related radial loading paths.

The third strategy represents the most suitable for a greater computational saving. As shown in Fig. 6, the three-dimensional strain space, defined by the orthogonal axes corresponding to the components $\bar{U}_{11}$, $\bar{U}_{12}$ and $\bar{U}_{22}$, is spanned by varying the time-like parameter $t$ for each radial direction, defined by angles $\theta$ (between 0° and 180°) and $\varphi$ (between 0° and 360°). The macroscopic right stretch tensor given as input at the RVE boundary value problem is defined by the following matrix:

$$\bar{U}(t) = \begin{bmatrix} 1+t\cos\varphi\sin\theta & t\cos\theta & 0 \\ t\cos\theta & 1+t\sin\varphi\sin\theta & 0 \\ 0 & 0 & 1 \end{bmatrix}. \qquad (20)$$

Each point corresponding to the restricted macroscopic first Piola-Kirchhoff stress is evaluated in a time step $t$ of the imposed radial loading path. The reduced database representing the nonlinear homogenized constitutive law can be hence extracted by imposing an appropriate number of radial loading paths and the time step ($\Delta t$). The database extraction procedure was implemented via a MATLAB script integrated with the finite element code COMSOL Multiphysics 5.4 by means of a parametric sweep able to vary $\theta$ and $\varphi$ through a prefixed range of values. Since the variational formulation of the adopted numerical environment is written in terms of the second Piola-Kirchhoff stress, $\bar{T}^{(2)} = \bar{F}^{-1}\bar{T}_R$. and of the Green-Lagrange strain, $\bar{E} = (\bar{F}^T \bar{F} - I)/2$, a database transformation step is needed before linking the reduced database to the noncritical subdomains.

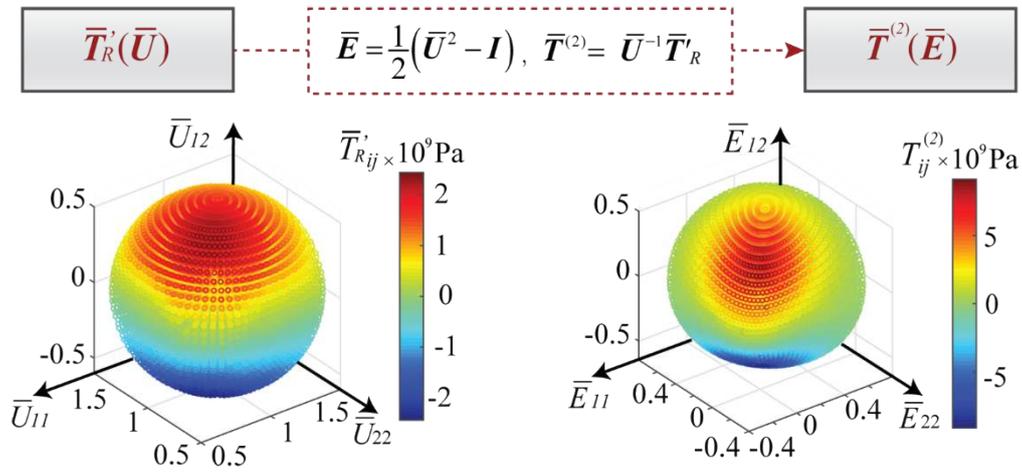

**Fig. 7** Example of transformation from the reduced database (on the left) to the transformed database (on the right), both depicted as point clouds in the related strain spaces.

For this reason, the previously discussed reduced database must be transformed using the following relations:

$$\bar{E} = \frac{1}{2}(\bar{U}^2 - I) , \qquad (21)$$

$$\bar{T}^{(2)} = \bar{U}^{-1}\bar{T}'_R , \qquad (22)$$

leading to a work-conjugated stress/strain database function of $\bar{E}$ and $\bar{T}^{(2)}$, as shown in Fig. 7. It worth noting that the fourth-order tangent constitutive tensor can be computed as follows:

$$\bar{C}^{(2)} = \frac{\partial \bar{T}^{(2)}}{\partial \bar{E}} . \qquad (23)$$

## 4 Numerical microscopic stability analysis

The above described multiscale approaches are here adopted in order to compute sequentially the principal solution path for the macroscopic models and the minimum eigenvalue of the microscopic structural stability functional, respectively, with the aim to perform the microscopic stability analysis of microstructured solids as described in Section 2. The two proposed multiscale strategies have been applied to differently arranged fiber-reinforced composite materials subjected to quasi-static and monotonically increasing macroscopic loads, in order to investigate their effectiveness in terms of both numerical accuracy and computational efficiency. Three numerical examples have been considered, the first one, reported in Section 4.1, being solved by means of the semi-concurrent multiscale approach described in Section 3.1, and the latter two, reported in Section 4.2, being analyzed via the hybrid hierarchical/concurrent multiscale strategy described in Section 3.2.

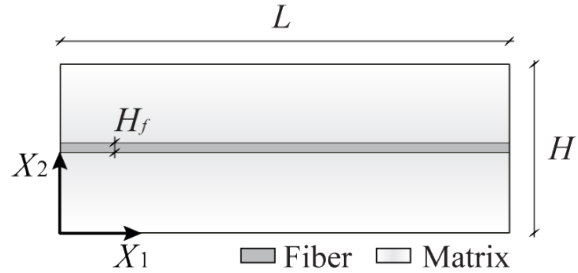

**Fig. 8** Micro-geometrical arrangement of the periodic RVE corresponding to a unidirectional fiber-reinforced composite material. The light gray area represents the soft material (matrix) and the dark gray area represents the stiff material (fiber).

### 4.1 Application of the semi-concurrent multiscale approach

The first multiscale application considers a composite material reinforced with continuous fibers arranged according to a unidirectional pattern. The associated layered microstructure is made of two sequentially repeated homogeneous layers, the thinner one representing the reinforcement and the other standing for the matrix, and it is characterized by the periodic unit cell sketched in Fig. 8. Moreover, the material interface between the two bulk phases is assumed to be perfect, meaning that no displacement discontinuity is allowed to occur. The constitutive law associated with the microscopic constituents is the neo-Hookean one and the corresponding strain energy density for plane strain deformations is of the following form:

$$W = \frac{\mu}{2}[F_{\alpha\beta}F_{\alpha\beta} - 2 - 2\ln J'] + \frac{k-\mu}{2}(J'-1)^2 \quad \alpha, \beta = 1, 2 \qquad (24)$$

where $J'$ is the determinant of the 2D deformation gradient tensor whose components are $F_{\alpha\beta}$, $\mu$ is the shear modulus of the solid at zero strain and the parameter $k$ plays the role of an equivalent 2D bulk modulus governing the material compressibility. The dimensions of the considered unit cell in the $X_1$ and $X_2$ directions (in the reference configuration), denoted as $L$ and $H$, respectively, are chosen such that their ratio $L/H$ is equal to 3 and $L$ is equal to 30 μm, whereas the thickness of the fiber (i.e. the reinforcement layer), indicated with $H_f$, is set as $0.025H$, associated with a fiber volume fraction of 2.5%. The dimensions of the macroscopic model in the $\bar{X}_1$ and $\bar{X}_2$ directions are $\bar{L}=240$ μm and $\bar{H}=40$ μm, respectively. The relative stiffness ratio $\mu_f/\mu_m$ between fiber and matrix is assumed to be equal to 200, with $\mu_m=807$ MPa, being the shear modulus at zero strain of the matrix material. The bulk modulus $k$ is assumed to be equal to $10\mu$ for both materials.

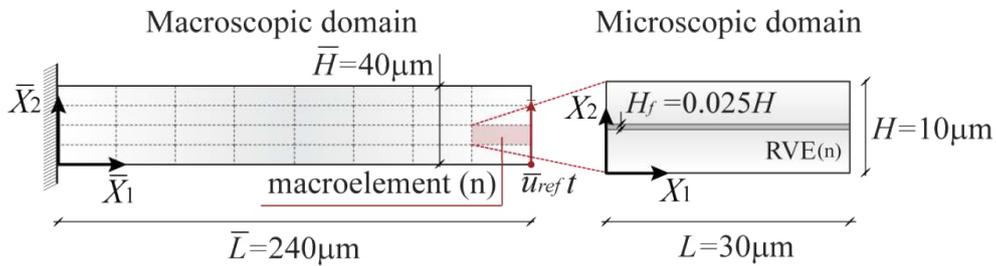

**Fig. 9** Cantilever composite beam reinforced with continuous fibers and subjected to a concentrated vertical force at the free end of the beam leading to a macroscopic displacement $\bar{u}_{ref}\,t$.

As shown in Fig. 9, the proposed semi-concurrent multiscale approach has been applied to model a cantilever composite beam reinforced with continuous fibers and subjected to a concentrated vertical force at the free end of the beam leading to a macroscopic displacement $\bar{u}_{ref}\,t$ (with $\bar{u}_{ref}=1$ μm), forced to grow monotonically with the time-like parameter $t$. A convergence analysis to the RVE size, reported in Fig. 10, has been incorporated in the numerical procedure to assess the local nature of the instability mode and, as a consequence, to verify the correspondence between the repeated unit cell (RUC) and the representative volume element (RVE).

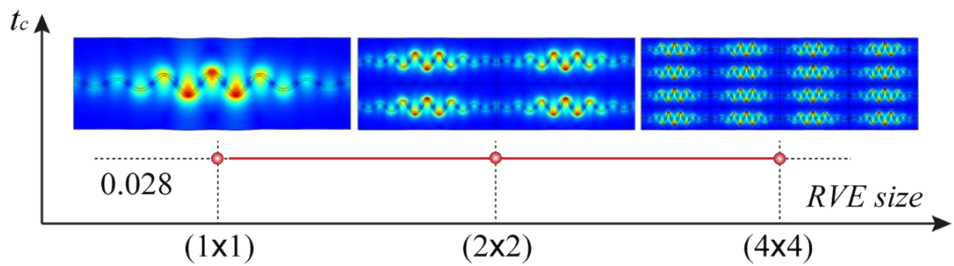

**Fig.10** Convergence analysis to the RVE size incorporated in the numerical procedure.

The microscopic instability critical load factor has been firstly evaluated by means of an extrapolation performed on the force-displacement curve and, successively, by means of a rigorous microscopic stability analysis following the theoretical development reported in Section 2. To evaluate the accuracy and the effectiveness of the adopted multiscale model, the obtained results have been compared with those coming from the stability analysis performed on a direct numerical model based on the explicit discretization of the heterogeneities of the composite microstructure.

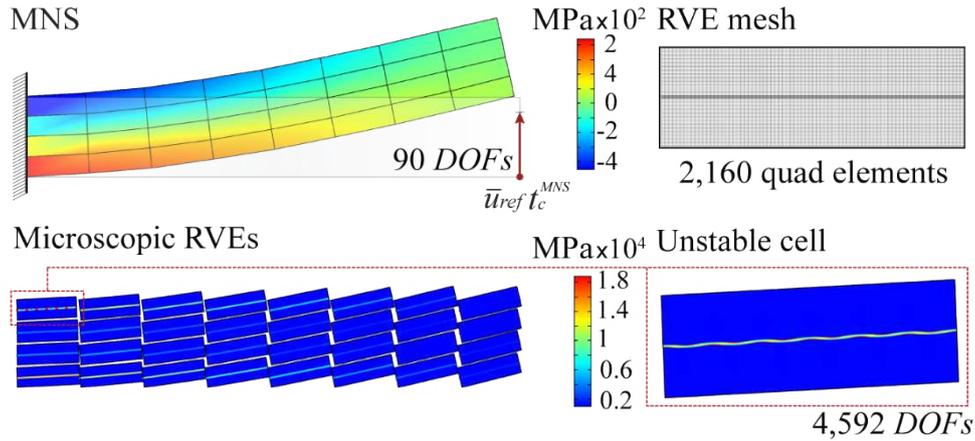

**Fig. 11** Multiscale numerical simulations (MNS) on a cantilever beam at the onset of microscopic instability corresponding to the critical load factor $t_c^{MNS}$.

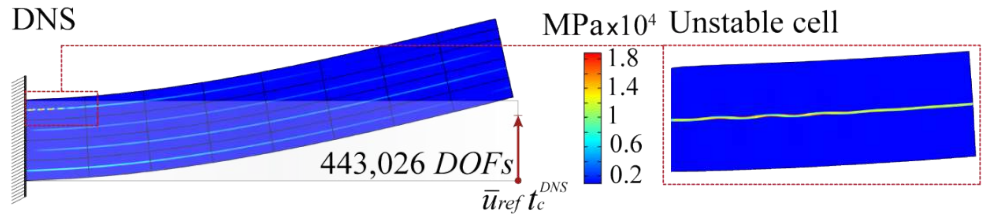

**Fig. 12** Direct numerical simulations (DNS) on a cantilever beam at the onset of microscopic instability corresponding to the critical load factor $t_c^{DNS}$.

In Fig. 11 the macroscopic multiscale model (MNS) is reported together with their linked microscopic representative volume elements (32 RVEs), giving also more information about the mesh discretization and the number of degrees of freedom involved (DOFs). Specifically, the adopted mesh of the macroscopic model involves 32 bilinear rectangular macro-elements and 90 degrees of freedom, while the microscopic RVEs involves 2,160 bilinear rectangular elements and 4,592 degrees of freedom. In the figure, the deformed configuration at the onset of microscopic instability induced by critical load factor $t_c^{MNS}$ and leading to a macroscopic vertical displacement $\bar{u}(t_c^{MNS}) = t_c^{MNS} \bar{u}_{ref}$ is also shown. In Fig. 12 the explicit model used for the direct numerical simulation (DNS) is reported, giving information about the number of degrees of freedom and showing the deformed configuration at the onset of the microscopic instability induced by critical load factor $t_c^{DNS}$ and leading to a macroscopic vertical displacement $\bar{u}(t_c^{DNS}) = t_c^{DNS} \bar{u}_{ref}$. Specifically, the adopted mesh involves 69,880 bilinear rectangular elements and 443,026 degrees of freedom. It is worth noting that both analyses predict the same cell (i.e. the upper left one) as the critical cell undergoing local instability.

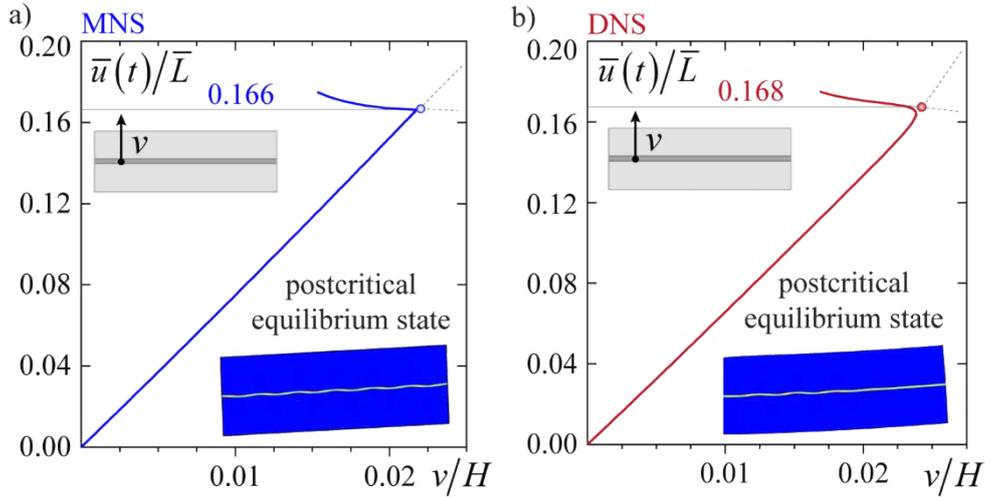

**Fig. 13** Comparison between the multiscale (a) and the direct (b) numerical simulations in terms of the extrapolated instability critical load factor for the cantilever beam.

The extrapolations of the critical load factors obtained by means of both numerical models are reported in Fig. 13. The force-displacement curves were obtained by plotting on the x-axis the local vertical displacement $v$ at the fiber-matrix interface of the most compressed unit cell involved by local instability, measured in the correspondence to the peak amplitude of the waveform, normalized with respect to the height of the unit cell $H$ ($v/H$) and by plotting on the y-axis the macroscopic vertical displacement $\bar{u}(t)$ normalized with respect to the macroscopic model length $\bar{L}$ ($\bar{u}(t)/\bar{L}$). With reference to Fig. 13a, the normalized critical load factor obtained by using the semi-concurrent multiscale method is equal to $\bar{u}(t_c^{MNS})/\bar{L} = 0.166$ leading to a critical load factor $t_c^{MNS} = 39.8$; while with reference to Fig. 13b, the normalized critical load factor obtained by using an explicit discretization of the heterogeneous microstructure is equal to $\bar{u}(t_c^{DNS})/\bar{L} = 0.168$ leading to an instability critical load factor $t_c^{MNS} = 40.3$. Comparing the results in terms of percent variation, being defined by the following relation:

$$e\% = \frac{t_c^{MNS} - t_c^{DNS}}{t_c^{DNS}}\%, \tag{25}$$

a relative error equal to -1.24% is obtained, thus confirming the good numerical accuracy of the adopted multiscale approach. This inaccuracy is essentially related to the use of homogenized moduli, computed by assuming periodic boundary conditions in every microscopic domain, even in regions where such conditions cease to hold due to the presence of boundary layer effects induced by constrained boundaries and concentrated loads.

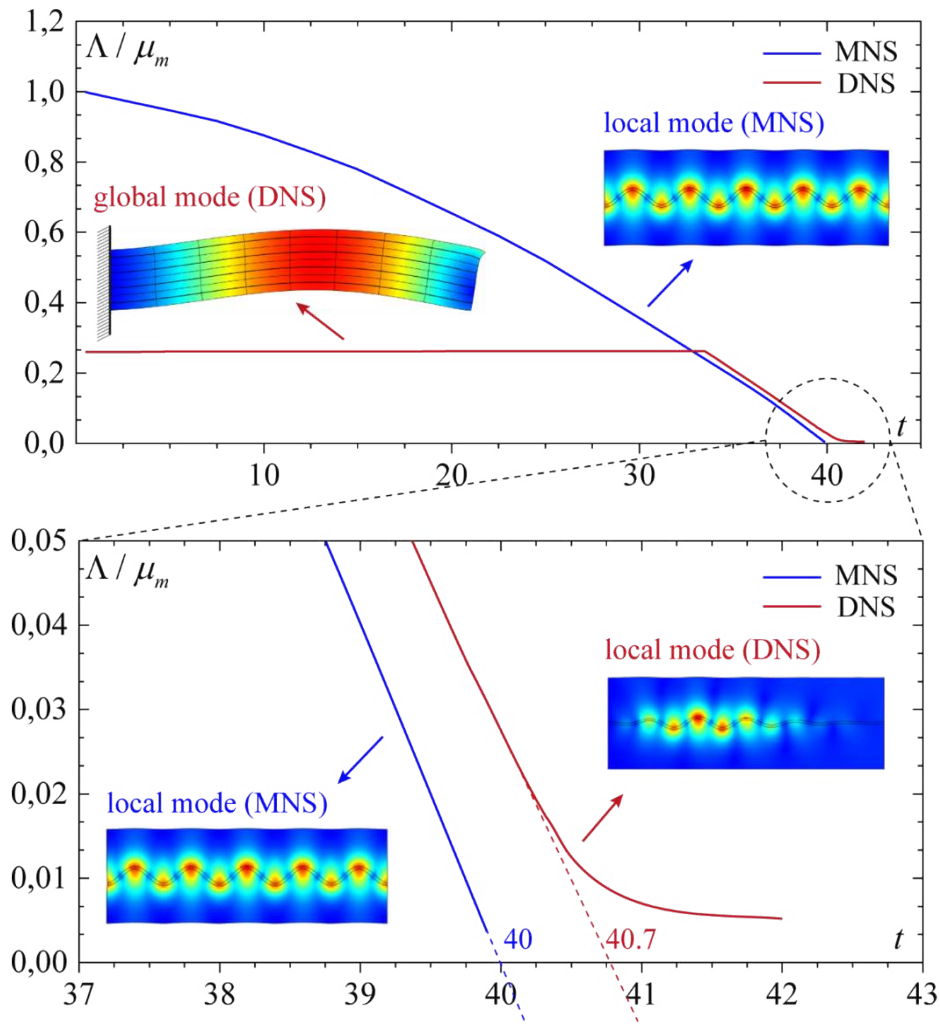

**Fig. 14** Normalized minimum eigenvalue plotted as a function of the time-like parameter *t* giving the critical load factor extracted by means of a rigorous microscopic instability analysis on the cantilever beam.

Subsequently, as shown in Fig. 14, the critical load factors have been evaluated by also performing a rigorous microscopic stability analysis, by which also the related mode shapes are obtained. The figure clearly shows that, with reference to the direct numerical analysis (red line), the primary instability mode is global in nature for lower values of the applied load, highlighting the competition between global and local instability modes for increasing load factors. On the contrary, with reference to the multiscale numerical simulation based on the semi-concurrent multiscale approach (blue line), the same behavior is not captured since, generally speaking, semi-concurrent approaches are not able to account for macroscopic boundary layer effect given by external constraints. In addition, as one can see in the zoomed view at the bottom of the figure, it worth noting that the minimum eigenvalue tends to zero without attaining this value, due to the presence of imperfections in the numerical model, so that a linear extrapolation technique was required to obtain this estimate. Specifically, respectively for the semi-concurrent multiscale and the direct numerical simulations, the following critical load factors have been estimated, $t_c^{MNS}=40$ and $t_c^{DNS}=40.7$. To summarize, resulting a relative percentage error less than 2% ($e\%=-1.72\%$), the semi-concurrent multiscale approach highlights a good prediction of the critical load factor in the case of microscopic instability, but the same accuracy is not exhibited in the prediction of the instability mode shapes because of the difficulty of semi-concurrent multiscale approaches to account for boundary layer effects.

**4.2 Application of the hybrid hierarchical/concurrent multiscale approach**

The two following hybrid hierarchical/concurrent multiscale applications consider a composite material reinforced with discontinuous fibers arranged according to a staggered pattern. The associated microstructure, characterized by the periodic unit cell sketched in Fig. 15, is made of two materials, the stiffer one representing the reinforcements in the form of elongated particles or short fibers and the softer one standing for the matrix.

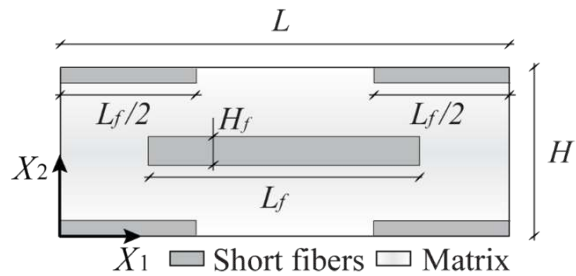

**Fig.15** Microgeometrical arrangement of the periodic RVE corresponding to a discontinuous fiber-reinforced composite material with a staggered pattern. The light gray area represents the soft material (matrix) and the dark gray area represent the stiff material (short fibers).

The constitutive law characterizing the mechanical behavior of the micro-constituents is the same of the first application reported in Section 4.1, see Eq.(24), together with the associated material parameters. The dimensions of the considered unit cell in the $X_1$ and $X_2$ directions are $H = (2 \cdot L_f \cdot H_f)/(L \cdot V_f)$ and $L = 400$ μm, the length of the short fibers is $L_f = 0.7L$ and their thickness is $H_f = L_f / 50$, leading to a fiber volume fraction $V_f$ of 12%. The dimensions of the macroscopic model in the $\bar{X}_1$ and $\bar{X}_2$ directions are $\bar{L} = 6400$ μm and $\bar{H} = 915$ μm, respectively. The relative stiffness ratio $\mu_f / \mu_m$ between fiber and matrix is assumed to be equal to 20, with $\mu_m = 807$ MPa, being the shear modulus at zero strain of the matrix material. The matrix bulk modulus is defined by $k_m = 2\mu_m$ and the fiber bulk modulus is defined by $k_f = 10\mu_f$.

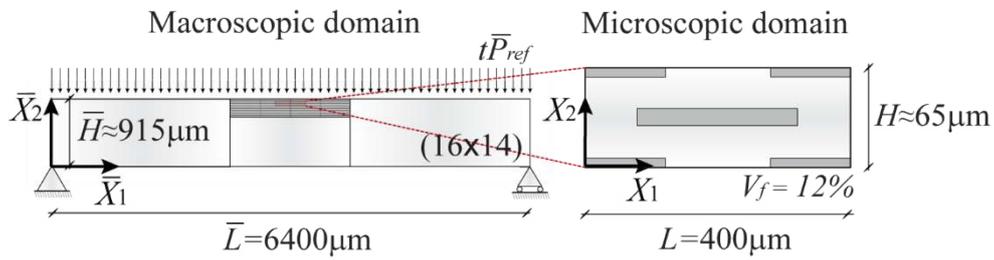

**Fig. 16** Simply supported composite beam reinforced with staggered discontinuous fibers and subjected to a distributed load $t\,\bar{P}_{ref}$.

In the first application, reported in Fig. 16, the proposed hybrid multiscale approach has been applied to analyze the structural behavior of a simply supported beam reinforced with staggered discontinuous fibers and subjected to a uniformly distributed load over the whole beam length equal to $t\bar{P}_{ref}$ with $\bar{P}_{ref} = 1\,e^7$ N/m. As in the previous application, the microscopic instability critical load factor has been firstly evaluated by means of an extrapolation performed on the force-displacement curve and, successively, by means of a rigorous microscopic stability analysis. Subsequently, these results have been validated by a suitable comparison with those coming from a direct numerical simulation.
The adopted multiscale model is reported in Fig. 17, which also indicates the number of degrees of freedom contained in the model, i.e. 80,564. In the figure, the deformed configuration at the onset of microscopic instability induced by the macroscopic distributed load $t_c^{MNS}\bar{P}_{ref}$ is also shown. In Fig. 18 the results obtained by the reference numerical simulation (DNS) are reported in terms of deformed configuration at the onset of the microscopic instability occurring at the load level $t_c^{DNS}\bar{P}_{ref}$. As expected, by considering an explicit discretization of the microstructure, a considerably higher number of degrees of freedom is involved, i.e. 1,134,136, leading to a huge computational cost.

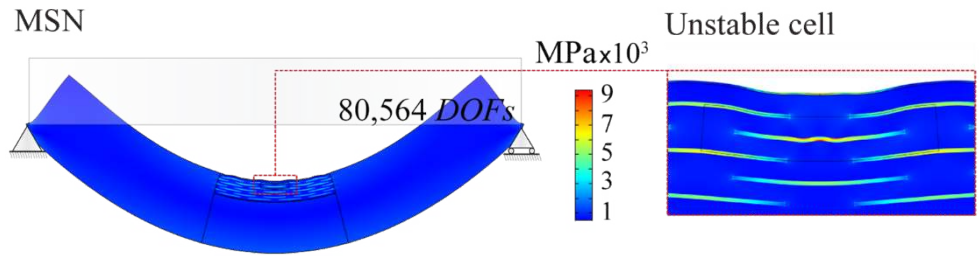

**Fig. 17** Multiscale numerical simulations (MNS) on the simply supported beam at the onset of microscopic instability corresponding to the critical load factor $t_c^{MNS}$.

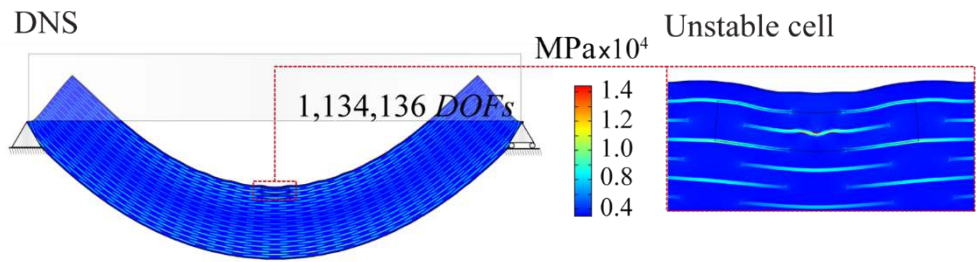

**Fig. 18** Direct numerical simulations (DNS) on the simply supported beam at the onset of microscopic instability corresponding to the critical load factor $t_c^{DNS}$.

The critical load factors are firstly determined by means of an extrapolation performed on the force-displacement curves obtained from both multiscale and direct numerical simulations, and reported in Fig. 19. The displacement $v^* = v - v_{mean}$ normalized with respect to the height of the unit cell $H$ ( $v^*/H$ ) is reported on the x-axis, whereas the load factor $t$ is plotted on the y-axis. Moreover, $v$ represents the vertical displacement at the central point of the upper fiber-matrix interface of the most compressed unit cell involved by local instability and $v_{mean}$ represents the mean vertical displacement of the midsection of the beam.

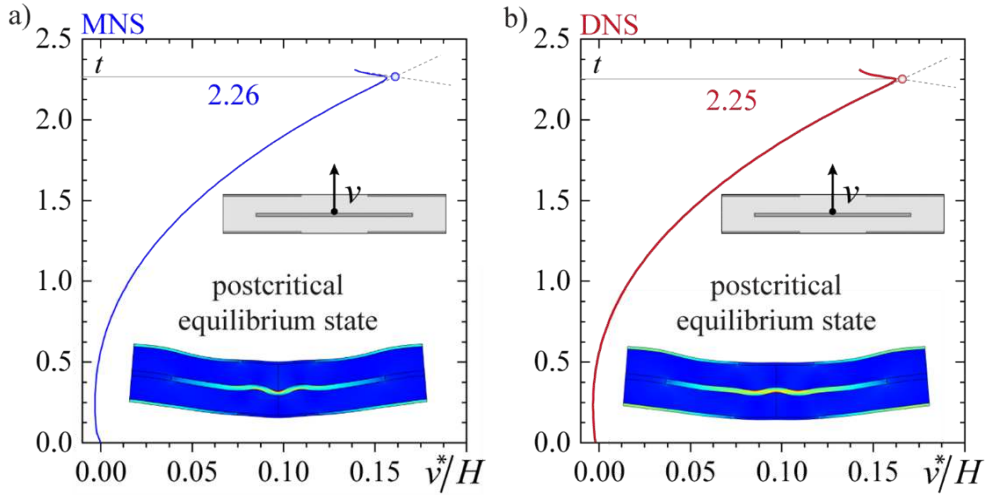

**Fig. 19** Comparison between the multiscale (a) and the direct (b) numerical simulations in terms of the extrapolation of the instability critical load factor for the simply supported beam.

With reference to Fig. 19a, the critical load factor obtained by using the hybrid multiscale method is equal to $t_c^{MNS}=2.26$; while with reference to Fig. 19b, the critical load factor obtained by using an explicit discretization of the heterogeneous microstructure is equal to $t_c^{DNS}=2.25$, leading to a relative percentage error equal to $e\%=0.44\%$. The excellent agreement between MNS and DNS results confirms the high accuracy of the proposed hybrid multiscale method in the determination of local instability critical load factors in discontinuous fiber-reinforced composite materials subjected to large deformations. To investigate more deeply the reliability of the proposed multiscale method, the critical load factors have been evaluated by also performing a rigorous microscopic stability analysis, whose results are reported in Fig. 20, also in terms of the related critical mode shapes.

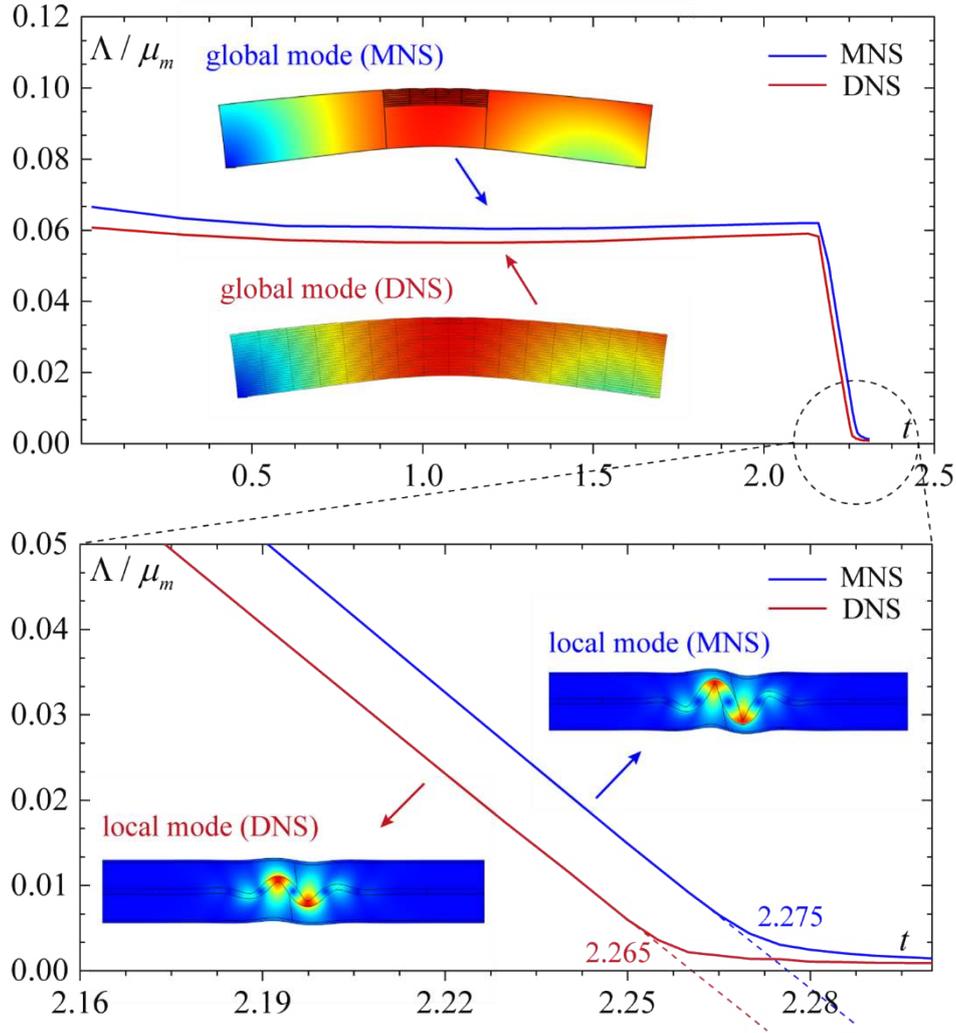

**Fig. 20** Normalized minimum eigenvalue plotted as a function of the time-like parameter *t* giving the critical load factor extracted by means of a rigorous microscopic instability analysis on the simply supported beam.

From this figure, it can be seen that, contrarily to the application reported in Section 4.1, both the direct numerical analysis (red line) and the multiscale numerical analysis (blue line) predict a primary instability mode of the global type for lower values of the applied load, also highlighting the competition between global and local instability modes for increasing load factors. Finally, a good agreement between the two analyses is reported in terms of shape of the global critical modes. Furthermore, it worth noting that the local mode shapes obtained via both analyses are in good agreement to each other and that the associated minimum eigenvalues tend to zero without attaining this value, due to the presence of imperfections in the numerical model. Thus, a linear extrapolation technique has been adopted allowing to determine the following critical load factor estimates: $t_c^{MNS} = 2.275$ and $t_c^{DNS} = 2.265$. On the basis of the very small resulting percentage error, equal to

$e\% = 0.44\%$, one can state that the hybrid multiscale approach highlights a very good prediction of the critical load factor in the case of microscopic instability and a very good prediction of the critical instability mode shape, showing a more pronounced accuracy in capturing the local instability modes than the global ones. These results demonstrate that the present hybrid multiscale method, contrarily to the previously applied semi-concurrent method, is able to account for boundary layer effects, thus leading to a more accurate prediction of the actual critical and post-critical behavior of composite materials.

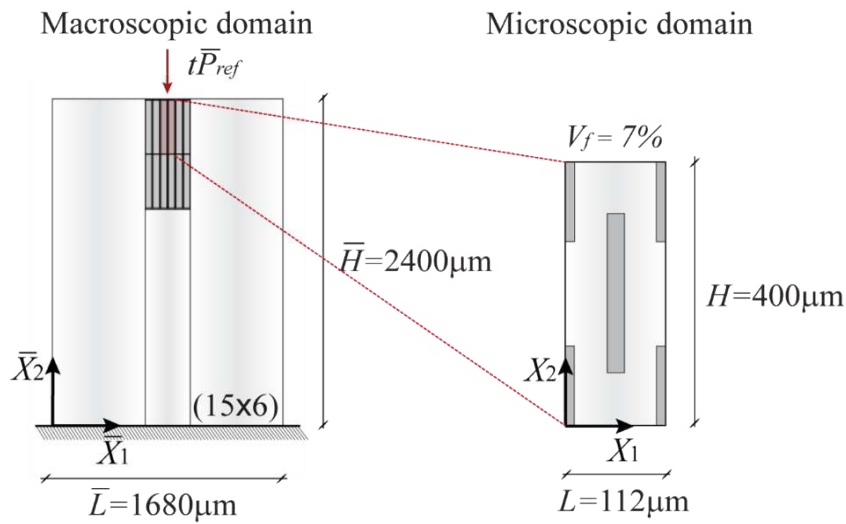

**Fig. 21** Clamped plate composite beam reinforced with staggered discontinuous fibers and subjected to a concentrated load $t\,\overline{P}_{ref}$ acting in the middle of the plate.

In the second application, reported in Fig.21, to further investigate the accuracy of the proposed hybrid multiscale approach under different loading conditions, a clamped plate reinforced with staggered discontinuous fibers and subjected to a concentrated load acting on the top boundary in correspondence of the midsection of the plate is considered. The dimension of the relevant unit cell in $X_1$ and $X_2$ directions are denoted as $L = (2 \cdot L_f \cdot H_f)/(L \cdot V_f)$ and $H = 400 \mu m$, where $L_f = 0.7L$ is the length of the short fibers, $H_f = L_f / 50$ is their thickness, whereas $V_f = 7\%$ represents the volume fraction of the stiffer phase. The dimensions of the macroscopic model in the $\overline{X}_1$ and $\overline{X}_2$ directions are $\overline{L} = 1680$ μm and $\overline{H} = 2400$ μm, respectively. The relative stiffness ratio $\mu_f / \mu_m$ between fiber and matrix is assumed to be equal to 10, with $\mu_m = 807$ MPa being the shear modulus at zero strain of the matrix material. The bulk modulus $k$ is assumed to be equal to 10 for both materials.

It worth noting that, generally speaking, if a material consisting of a thin stiff layer and a softer substrate is subjected to a sufficiently large compressive load, a buckling or wrinkling surface instability can occur, as shown in [52]. Because of the interest in the microscopic instability of the short fibers embedded in a soft matrix, of the occurrence of any type of surface instability is prevented by inserting a stiffer thin layer on the top side of the plate over its entire width. The thickness of this layer is equal to $H_f/2$, and its material is characterized by a linear elastic constitutive response with elastic modulus equal to 200GPa and $\upsilon=0$.

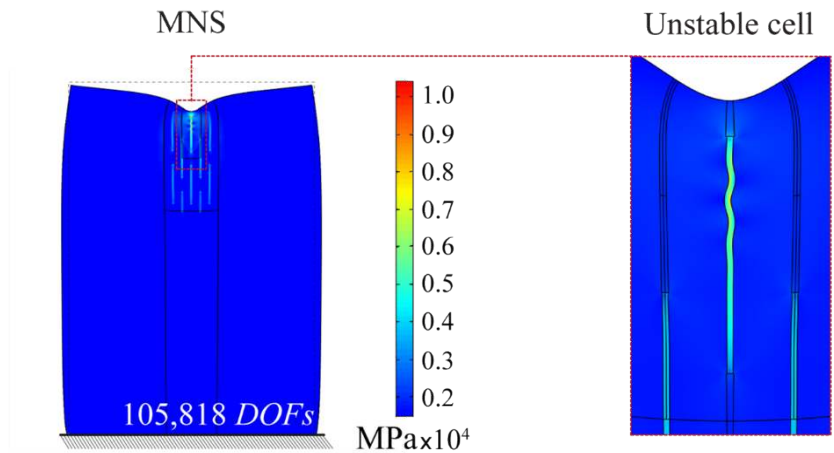

**Fig. 22** Multiscale numerical simulation (MNS) on the clamped plate at the onset of microscopic instability corresponding to the critical load factor $t_c^{MNS}$.

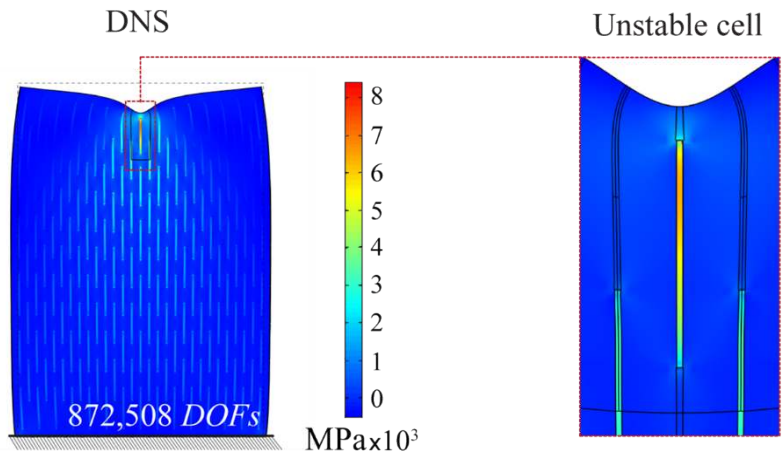

**Fig. 23** Direct numerical simulation (DNS) on the clamped plate at the onset of microscopic instability corresponding to the critical load factor $t_c^{DNS}$.

The deformed configurations corresponding to the onset of microscopic instability obtained by means of a multiscale numerical simulation (MNS) and by means of a direct numerical simulation (DNS) are respectively reported in Fig. 22 and Fig. 23. Also in this application, it is clearly highlighted that the adopted hybrid multiscale approach leads to a numerical model with a lower number of DOFs and, as a consequence, it leads to a sensibly lower computational effort with respect to the direct approach. Specifically, the multiscale finite element model is characterized by 105,818 DOFs, whereas the explicit finite element model possesses more than 8 times the number of multiscale DOFs.

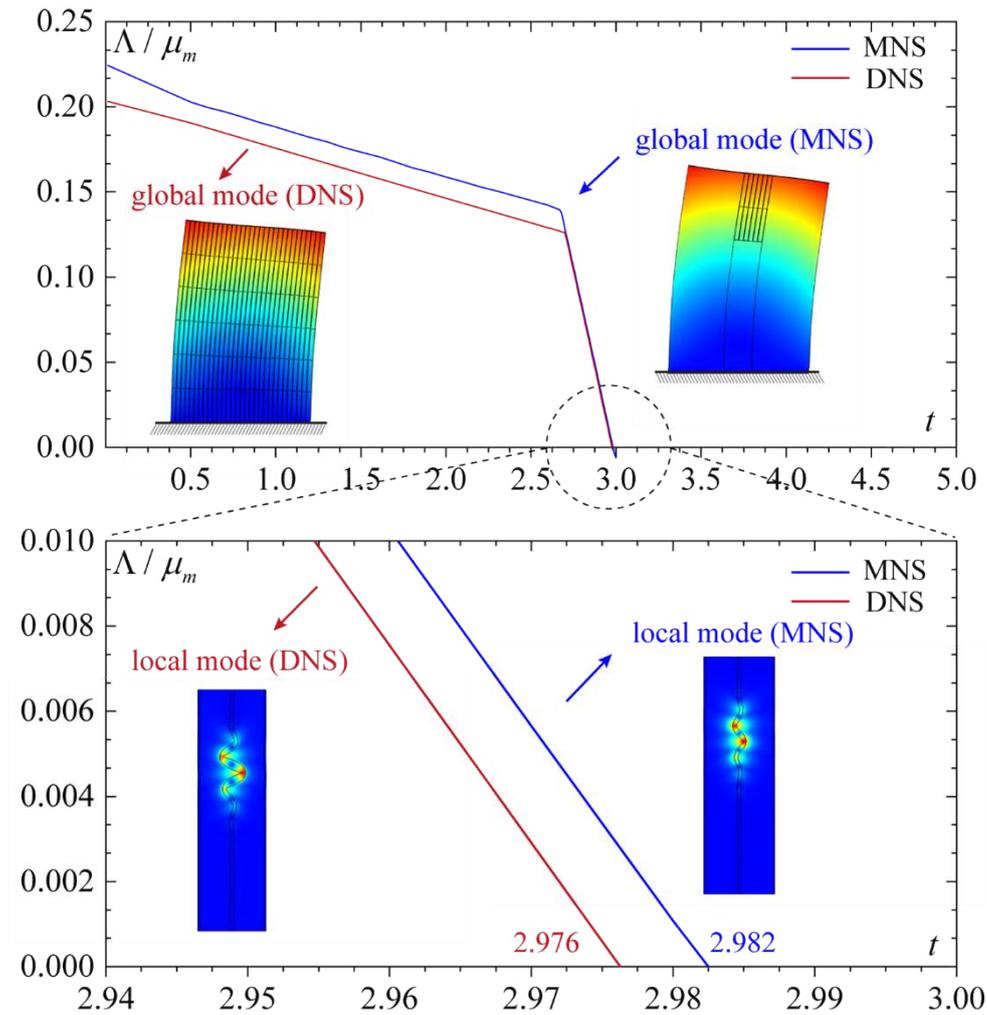

**Fig. 24** Normalized minimum eigenvalue plotted as a function of the time-like parameter *t* giving the critical load factors extracted by means of a rigorous microscopic instability analysis on the clamped plate.

The microscopic instability critical load factors have been evaluated by means of a microscopic stability analysis, whose results are reported in Fig. 24, also in terms of the related critical mode shapes. The following critical load factors are obtained via the multiscale numerical simulation

(MSN) and the direct numerical simulation (DNS), respectively i.e. $t_c^{MNS} = 2.976$ and $t_c^{DNS} = 2.982$, leading to a relative percentage error equal to $e\% = 0.2\%$. Also in this application, the critical mode shapes, obtained by means of MNS and DNS approaches, are perfectly coincident, thus allowing to further assess the good accuracy of the proposed hybrid multiscale method in the prediction of the either local or global critical mode shapes.

Definitely, the obtained results are in good agreement with that obtained above in the first hybrid multiscale application, leading to the consideration that the proposed hybrid multiscale approach has been shown more effective than the semi-concurrent approach for obtaining accurate estimates of the critical load levels associated with microscopic instabilities as well as of the related critical mode shapes, being able to effectively account for the boundary layer effect and, as a consequence, for the competition between global and local critical mode shapes during the loading phase.

## 5 Conclusions

This work is a contribution to the analysis of instabilities in finitely strained composite materials whose micro-constituents are characterized by nonlinear constitutive relations. Specifically, their macroscopic mechanical response was investigated by using nonlinear homogenization techniques and advanced computational multiscale strategies. Such advanced modeling techniques constitute an effective tool to analyze the structural response of a large class of composite materials with complex microstructures, providing a link between the macroscopic behavior and the underlying microstructural phenomena, accounting for different types of nonlinearity, such as large deformation, fracture, contact or instability (or also a combination of these).

Hence, we proposed two multiscale modeling strategies to investigate the microstructural instability in locally periodic fiber-reinforced composite materials subjected to general loading conditions in a large deformation context. The first adopted strategy is a semi-concurrent multiscale method consisting in the on-the-fly derivation of the macroscopic constitutive response of the composite structure (as in the most of FE$^2$-like methods), here used in conjunction with a microscopic stability analysis and integrated within a new two-way computational homogenization scheme. The second approach is a novel hybrid hierarchical/concurrent multiscale approach able to combine the advantages inherent in the use of hierarchical and concurrent approaches, and based on a two-level domain decomposition solution scheme. The main goal of the proposed multiscale strategy is to combine a hierarchical multiscale approach, in domains wherein the assumption of scale separation is fulfilled (homogenized domains), with a concurrent approach, in remaining domains for which this condition is no longer fulfilled (fine-scale domains) due to the presence of high macroscopic stress or strain gradients. Specifically, fine-scale domains are modeled with an explicit description of all its microstructural details. Homogenized domains, instead, use information extracted by the solution of a microscopic BVP defined over an RVE; in particular, a microscopically informed macroscopic constitutive relation in the form of a macro-stress/macro-strain database has been implemented in the finite element model in conjunction with an interpolation method.

The viability and accuracy of the proposed multiscale approaches in the context of the microscopic stability analysis of composite materials have been appropriately evaluated through comparisons with reference direct numerical simulations.

As a first numerical application, the proposed semi-concurrent multiscale approach has been used to predict local instabilities in a cantilever composite beam reinforced with periodically distributed continuous fibers. In particular, an estimate of the instability critical load has been initially determined via an extrapolation performed on the numerically obtained load-displacement curve. The numerical accuracy of the multiscale approach has been assessed by the comparison with a direct numerical simulation, based on a fully meshed model. The error on the estimated critical load between the two analyses is of about 1%, mainly due to both macroscopic gradients in stresses and strains and boundary layer effects. Then, the critical load factor has been also estimated by rigorous microscopic stability analysis. The good numerical accuracy of the semi-concurrent method is confirmed by a small error on the critical load, less than 2%. Interestingly, with reference to the direct numerical analysis, the competition between global and local modes has been highlighted; on the contrary, with reference to the semi-concurrent multiscale approach, the same behavior is not captured since semi-concurrent approaches are not capable to evaluate boundary layer effects given by external constraints. In addition, the predicted local mode is not perfectly coincident with the reference one, due again to the fact that boundary layer effects are not accounted for.

As a second numerical application, the novel hybrid multiscale approach has been applied to two different composite structures reinforced with periodic staggered discontinuous fibers. The first case concerns the bending problem for a simply supported composite beam. The numerical accuracy of the proposed multiscale approach has been assessed by the comparison with a direct numerical simulation. A very small error on the estimated critical load between the two analyses is found, less than 0.5%. The critical load factor has been also estimated by a rigorous microscopic stability analysis. In this case, the competition between global and local modes is captured by both direct and multiscale analyses and the excellent numerical accuracy of the hybrid hierarchical/concurrent multiscale method is confirmed by the very small error on the extrapolated critical load, of about 0.5%. Moreover, the predicted local mode is perfectly coincident with the reference one, due to the fact that boundary layer effects are accounted for. The last case, concerns a clamped composite plate under a concentrated load. The great numerical accuracy of the hybrid multiscale approach is confirmed by the very small error on the critical load, of about 0.2%. Moreover, also in this case, the predicted local mode is perfectly coincident with the reference one. In conclusions, suitable comparisons with direct numerical simulations have shown that the hybrid multiscale approach is more accurate in obtaining both critical load levels and related microscopic instability modes, by virtue of its capability of capturing both macroscopic gradients in stresses and strains and boundary layer effects. In addition, from the point of view of computational effort, the proposed hybrid multiscale model is, in general, more efficient compared with both direct and semi-concurrent models, but the claimed efficiency is strictly related to the a-priori identification of critical and noncritical domains. However, the prior determination of regions undergoing local instabilities is not trivial in the presence of more complex geometries and/or loading configurations. Therefore, a further extension of this work could concern the development of a more general hybrid hierarchical/concurrent multiscale method equipped with adaptive capabilities, able to deal with microscopic instability phenomena in a more efficient manner. Specifically, this envisaged adaptive hybrid multiscale approach, inspired from similar approaches already developed in a small deformation framework [17,53–55], could be developed starting from the definition of suitable zooming-in criteria to adaptively trigger the identification and the subsequent insertion of new critical regions, with the explicit aim of injecting a detailed microstructure exclusively in those zones characterized by an incipient failure state.

**Acknowledgements**

Fabrizio Greco and Paolo Lonetti gratefully acknowledge financial support from the Italian Ministry of Education, University and Research (MIUR) under the P.R.I.N. 2017 National Grant "Multiscale Innovative Materials and Structures" (Project Code 2017J4EAYB; University of Calabria Research Unit).